%% file: main.tex
\documentclass[onecolumn,11pt,a4paper,reqno]{article}
\input{settings.tex}

\begin{document}
\input{preamble.tex}

\begin{abstract} 
The Demazure product, also called the 0-Hecke product, is an associative operation on Coxeter groups with interesting properties and applications. In \cite{demAhopping} it was shown that the Demazure product of two permutations can be described purely combinatorially: using only their one-line notation and not relying on reduced words. In this paper, we extend this to Coxeter groups of type $D$.
\end{abstract}

\section{Introduction}
Coxeter groups play a significant role in both the representation theory of Lie groups and the geometric study of their associated flag and Schubert varieties. An intriguing associative operation known as \newword{Demazure product} \cite{Dem94} (also known as the $0$-\newword{Hecke} or \newword{greedy product}) is defined on Coxeter groups.

Let $W$ be a Coxeter group with generating set $S$. It has relations of the the form 
\begin{equation}\label{eqn:Coxeter_braid}
(st)^{m_{st}}=id, \quad s,t\in S
\end{equation}
for some values $m_{st}\in \Z_{>0}\cup\{\infty\}$ where $m_{st}=1$ if and only if $s=t$.
Coxeter groups are associated with a length function $\ell:W\rightarrow \Z_{\geq 0}$ and a partial order $\leq$, known as the Bruhat order. For in-depth information on the properties of the Coxeter group, refer to \cite{Bjorner-Brenti05}. The \newword{Coxeter monoid} (also termed $0$ -\newword{Iwahari-Hecke monoid}) is defined on $W$ as the monoid formed by generators $S$ using a product $\star$ that obeys the Coxeter braid relations provided in Equation \eqref{eqn:Coxeter_braid} for $s\neq t$, along with the rule $s\star s=s$ for each $s\in S$. The product $\star$ is known as the Demazure product, initially investigated by Norton within the realm of Hecke algebras in \cite{No79}. As sets, we have $W=\langle S,\star\rangle$. An expression $w=s_1\dots s_k$ is termed \newword{reduced} if $\ell(w)=k$, and in this case $w$ cannot be formed using fewer than $k$ generators from $S$. We say $u\leq w$ in \newword{Bruhat order} if there exists a reduced expression for $u$ that is a subsequence of a reduced expression for $w$ (\cite{Bjorner-Brenti05}, Theorem 2.2.2). The following lemma notes some elementary facts about the Coxeter monoid.

\begin{lem}\cite[Lemma 1.3, Corollary 1.4]{No79}
\label{lemma:moniod_properties}
Let $W$ be a Coxeter group with generating set $S$.  Then the following are true:
\begin{enumerate}
    \item Let $(s_1,\dots,s_k)$ be a sequence of generators in $S$.  Then
    $$s_1\dots s_k\leq s_1\star\dots \star s_k$$ with equality if and only if $(s_1,\dots,s_k)$ is a reduced expression.
    \item For any $s\in S$ and $w\in W$,
$$s\star w=\begin{cases} w & \text{if}\quad \ell(sw)< \ell(w)\\ sw & \text{if}\quad \ell(sw)>\ell(w).\end{cases}$$
\end{enumerate}
\end{lem}

 The \newword{Bruhat interval} is defined as $[u,w]:=\{v\in W \ |\ u\leq v\leq w\}$. In cases where $u=id$, the interval $[id,w]$ is referred to as \newword{lower interval} of $w$. There exists an elegant interpretation of $w \star v$ for any $w,v\in W$ in terms of Bruhat intervals:

\begin{prop}\cite[Lemma 1]{He09}, \cite[Proposition 8]{Kenny14}\label{prop:interval_product}
For any $w,v\in W$, the lower interval $$[id,w\star v]=\{ab \ |\ a\in [id,w], b\in[id,v]\}.$$
\end{prop}

This product has found extensive application and has been studied in multiple fields that use Coxeter groups \cite{gay:hal-01691266}, \cite{KNUTSON2004161}, \cite{Larson19}, \cite{OR}, \cite{Pflueger}, \cite{Richardson-Springer90}. In \cite[Fact 2.4]{Chan-Pflueger23}, Chan and Pflueger describe how the Demazure product of permutations can be expressed using rank functions. Subsequently, in \cite{Pflueger}, Pflueger demonstrates that this description extends to broader categories, specifically almost-sign-preserving permutations. Meanwhile, Tiskin in \cite{Tiskin15} provides a recursive method to compute the Demazure product using the Monge matrix (rank matrix) of a permutation and using tropical operations.


In \cite{demAhopping} a method to calculate this product via a new tool termed \newword{hopping operator} was introduced for Coxeter groups of type $A$ and type $B$. This operator enables us to compute the Demazure product without resorting to reduced words, lengths, or simple transpositions but just using the one-line (window) notation. In this paper, we will extend this result to type $D$ Coxeter groups. In Section 2 we review the hopping operator for types $A$ and $B$. In Section 3 we go over some properties and tools for the type $D$ Coxeter groups. In Section 4 we present and prove our main result.

\section{Hopping operator for type $A$ and $B$}
In this section we review the results of \cite{demAhopping}. The examples are also borrowed from \cite{demAhopping}. 

\subsection{Hopping in type $A$.}

The permutation group $S_n$ is a Coxeter group of type $A$ with simple generating set $S=\{s_1,\dots,s_{n-1}\}$ satisfying the relations $s_i^2=id$ and 
\begin{equation}\label{eqn:type_A_braid}
(s_is_j)^2=id\ \text{ if $|i-j|>1$\ and}\ (s_is_{i+1})^3=id.
\end{equation}  
The generator $s_i$ corresponds to the simple transposition $(i,i+1)$ and, for $w\in S_n$, let $w=w(1)w(2)\dots w(n)$ denote the permutation in one-line notation.  We denote the set of integers $\dbrack{n}:=\{1,2,\dots, n\}$ and use the notation $L:=[a_1,\dots,a_k]$ to denote an ordered subset of $\dbrack{n}$ (without repetition). Beware that within this notation, $[a]\neq \dbrack{a}$.  

\begin{defn}\label{def:hopping_typeA}
For $t \in \dbrack{n}$ and $L$ an ordered subset of $\dbrack{n}$, the \newword{hopping operator} $$h_{t,L}:S_n\rightarrow S_n$$ acts on a permutation $w$ according to the following algorithm: Scan to the right (within the one-line notation of $w$) of $t$ and look for the element furthest to the right in $L$ that is greater than $t$. If it exists, swap $t$ and that element, replace $w$ with the resulting permutation, and repeat. The algorithm ends when there are no elements of $L$ within $w$ to the right of $t$.
\end{defn}

For example, take $n=8$ and $w = 891726435$. Then $h_{1,[2,3,4,5,6,7,8]}(w) = 897625431$ is obtained by the following process:
$$89\textbf{17}26435 \rightarrow 897\textbf{1}2\textbf{6}435 \rightarrow 89762\textbf{1}43\textbf{5} \rightarrow 897625431.$$

Given $w \in S_n$, we let $w \nwarrow a$ stand for the subword of the one-line notation of $w$ obtained by restricting to the subword strictly left of $a$, then removing all elements smaller than $a$.

\begin{thm}[\cite{demAhopping}, Theorem 12]\label{thm:Amain}
For any $w,v\in S_n$, we have $w \star v = h_{n-1,w{\nwarrow n-1}} \dots h_{2,w{\nwarrow 2}} h_{1,w{\nwarrow 1}}(wv)$.
\end{thm}

For example, let $w=6541723$ and $v=5436217$.  The usual product of these two permutations is $wv=7142563$.  The Demazure product $w\star v$ corresponds to the sequence of hopping operators
$$h_{5,[6]}h_{4,[6,5]}h_{3,[6,5,4,7]}h_{2,[6,5,4,7]}h_{1,[6,5,4]}$$
acting on the usual product $wv$.  Applying each of the hopping operators in order to $wv$, we get 
\begin{align*}
7142563 \xrightarrow[h_{1,[6,5,4]}]{}  7452613 \xrightarrow[h_{2,[6,5,4,7]}]{} & 7456213\\  
& \xrightarrow[h_{3,[6,5,4,7]}]{}  7456213  \xrightarrow[h_{4,[6,5]}]{}  7564213 \xrightarrow[h_{5,[6]}]{} 7654213.
\end{align*}

\subsection{Hopping in type $B$.}
\label{subsec:typeB}

The group of signed permutations is also known as Coxeter groups of type $B$ (or equivalently, type $C$).  Signed permutations can be viewed as a permutation subgroup of $S_{2n}$. Let $\{s'_1,\dots,s'_{2n-1}\}$ be the simple generators of the permutation group $S_{2n}$.  We define $B_n$ to be the subgroup of $S_{2n}$ generated by $S:=\{s_1,\dots,s_n\}$ where
\begin{equation}\label{eqn:BtoA}
s_i:=s'_i\, s'_{2n-i}\ \text{for $1\leq i< n$ and}\quad s_n:=s'_n.
\end{equation}

We use the convention of \cite{billey2000singular} when working with simple transpositions of type $B$. As a Coxeter group, the generators $s_1,\dots,s_{n-1}$ of $B_n$ satisfy the same relations as in type $A$ (see Equation \eqref{eqn:type_A_braid}) with the last generator $s_n$ satisfying: 
$$(s_is_n)^2=id\ \text{ for $1\leq i<n-1$ and } (s_{n-1}s_n)^4=id.$$

The one-line notation for type $B$ works as follows. The convention we use here is slightly different from that of Bjorner and Brenti in \cite{Bjorner-Brenti05} in the sense that $s_n$ plays the role of $s_0$.

\begin{defn}
A \newword{signed permutation} of type $B_n$ is a permutation of $\dbrack{n}$ along with a sign of $+$ or $-$ attached to each number.
\end{defn}

Let $\pm \dbrack{n}$ denote the set $\dbrack{n} \cup -\dbrack{n}$, where $-\dbrack{n} := \{-1,\dots,-n\}$.   We impose the total ordering on $\pm\dbrack{n}$ given by:
$$1 < 2<\dots < n < -n <  \dots<-2 < -1.$$ 

By \newword{unfolding} of a signed permutation, we mean the following: to the right of $w$, attach a reverse ordered copy of $w$ with the signs flipped to get a permutation of $\pm \dbrack{n}$. For example the unfolding of $[4,-2,3,-1]$ is $$[4,-2,3,-1,1,-3,2,-4]$$ and the corresponding permutation of $\dbrack{8}$ is $[4,7,3,8,1,6,2,5].$

\begin{defn}[\cite{demAhopping}]
Let $t \in \pm \dbrack{n}$ and $L$ be an ordered subset of $\pm \dbrack{n}$ (without repetition).  The \newword{hopping operator} 
 acts on a signed permutation $w$ by the following algorithm: scan to the right (within the unfolding of $w$) of $t$ and look for the element furthest to the right in $L$ that is greater than $t$. If it exists, say $q$, then swap $t$ and $q$ and also swap $-t$ with $-q$ (unless $t = -q$).  Replace $w$ with the resulting unfolded signed permutation and repeat. The algorithm ends when there are no elements of $L$ within $w$ to the right of $t$.
\end{defn}

For example, let $w = [-4,5,3,-1,-2,-6]$ and we will calculate how the hopping operator $h_{1,[-3,4,-5, 6]}$ acts on it. First we unfold $w$, which gives $$\unfold(w) = [-4, 5, 3, -1,-2,-6, 6, 2, 1, -3, -5, 4].$$ To the right of $1$ in $\unfold(w)$ we have $[-3,-5,4]$ among the elements of $L$. We first swap $1$ for $-5$ since $-5$ is the rightmost element of $L$ that exists here (we also swap $-1$ and $5$). This gives us $$[-4, -1, 3, 5,-2,-6, 6, 2, -5, -3, 1, 4].$$ After that, we again scan to the right of $1$ to find $[4].$ Then we swap $1$ with $4$, to get $$[-1, -4, 3, 5,-2,-6, 6, 2, -5, -3, 4, 1].$$ So the signed permutation we end up with is $[-1,-4,3,5,-2,-6]$.

For any $w\in B_n$ and $i>0$, define $w \nwarrow i$ as the subword of $\unfold(w)$ obtained by restricting to the numbers to the left of $i$ contained in $(i,-i]$.

\begin{thm}[\cite{demAhopping}, Theorem 20]\label{thm:Bmain}
For any $w,v\in B_n$, we have $w \star v = h_{n-1,w{\nwarrow n-1}} h_{2,w{\nwarrow 2}} h_{1,w{\nwarrow 1}}(wv)$.
\end{thm}

For example, let $w=[-5,3,1,-2,4]$ and $v=[-4,2,-1,-3,5]$ in $B_5$. The unfolding of $w$ is $[-5,3,1,-2,4,-4,2,-1,-3,5]$. We start with the usual product $wv=[2, 3, 5, -1, 4].$  We apply the sequence of hopping operators 
$$h_{5,[-5]}h_{4,[-5]}h_{3,[-5]}h_{2,[-5,3,-2,4,-4]}h_{1,[-5,3]}$$
to $wv$ giving:

\begin{align*}
& [2,3,5,-1,4,-4,1,-5,-3,-2] \xrightarrow[h_{1,[-5,3]}]{} [2,3,-1,5,4,-4,-5,1,-3,-2]  \xrightarrow[h_{2,[-5,3,-2,4,-4]}]{} \\
& [-2,3,-1,5,-4,4,-5,1,-3,2] \xrightarrow[h_{3,[-5]}]{} [-2,-5,-1,-3,-4,4,3,1,5,2]  \xrightarrow[h_{4,[-5]}]{} \\
& [-2,-5,-1,-3,-4,4,3,1,5,2] \xrightarrow[h_{5,[-5]}]{} [-2,-5,-1,-3,-4,4,3,1,5,2],
\end{align*}

then \cref{thm:Bmain} implies that $w\star v=[-2, -5, -1, -3, -4]$.

\begin{rmk}
It is natural to ask why not use the exact same proof techniques to try to extend the type $B$ result to type $D$? It turns out that unfolding behaves drastically different for type $B$ and type $D$, which serves as an obstacle for techniques to carry over. We go over this in more detail in \cref{subsec:unfolding}.
\end{rmk}

\section{Going over various properties of type $D$ elements}

In this section, we will develop various tools that we will use for proving our main result in the following section. The Coxeter diagram for type $D$ is given in \cref{fig:dn}.

\input{figures/typeDdiagram}

Notice from the figure that we are using the same convention as in \cite{billey2000singular} when working with Coxeter groups of type $D$. The set of simple transpositions is $S := \{s_1,\dots,s_n\}$ where $s_1,\dots,s_{n-1}$ satisfy the same relations as in type $A$.

\subsection{One-line notation for type $D$}
Similar to the type $A$ and $B$ cases, the elements of the Coxeter group $D_n$ can be expressed using a one-line notation as the definition below. Similarly to how the case of type $B$ was treated in \cite{demAhopping}, let $\pm \dbrack{n}$ denote the set $\dbrack{n} \cup -\dbrack{n}$, where $-\dbrack{n} := \{-1,\dots,-n\}$.   We impose the total ordering on $\pm\dbrack{n}$ given by:
$$1 < 2<\dots < n < -n <  \dots<-2 < -1.$$ 

\begin{defn}
A \newword{signed permutation} is a permutation of $\dbrack{n}$ along with a sign of $+$ or $-$ attached to each number. A signed permutation is called \newword{even} if it has exactly even number of negative signs.
\end{defn}

For notational simplicity, we drop the $+$ signs from the one-line notation. The convention we use here is slightly different from that of Bjorner and Brenti in \cite{Bjorner-Brenti05} in the sense that $s_n$ plays the role of $s_0$. The product structure on (even) signed permutations is simply the usual composition of permutations with the added condition that $w(-i)=-w(i)$.

Here is how the simple transpositions act on even-signed permutations, which allows us to go between reduced words of elements of Coxeter group $D_n$ and one-line notation (even-signed permutations).

\begin{itemize}
    \item $s_i w$ for $i < n$ swaps the entries $i$ and $i+1$ (the signs remain at their positions).
    \item $w s_i$ for $i < n$ swaps the entries at position $i$ and $i+1$ (the signs stay with the entries they were paired with).
    \item $s_n w$ swaps the entries $n-1$ and $n$, then flips the sign of $n-1$ and the sign of $n$.
    \item $w s_n$ swaps the positions $n-1$ and $n$, then flips the sign at position $n-1$ and position $n$.
\end{itemize}

For example, for $n=5$, we have the following:
\begin{itemize}
    \item $s_3 * [1,2,-3,4,-5] = [1,2,-4,3,-5]$,
    \item $[1,2,-3,4,-5] * s_3 = [1,2,4,-3,-5]$,
    \item $s_5 * [1,2,-3,4,-5] = [1,2,-3,-5,4]$,
    \item $[1,2,-3,4,-5] * s_5 = [1,2,-3,5,-4]$.
\end{itemize}

Consider $w = s_1s_2s_1s_3s_5s_3s_2$. To get its one-line notation, we do $w = s_1s_2s_1s_3s_5s_3s_2 * [1,2,3,4,5] = s_1s_2s_1s_3s_5s_3*[1,3,2,4,5] = \dots = [3,-5,2,4,-1]$.


\subsection{Combinatorial unfolding of the even signed permutations}
\label{subsec:unfolding}

A (simply-laced) Coxeter–Dynkin diagram that has a symmetry can be quotiented by the symmetry, resulting in a new, potentially multi-laced diagram, via the process called 
\newword{geometric folding} \cite{Zuber1998}. This gives a natural way to embed signed permutations in $S_{2n}$. But for even signed permutations, geometric folding does not give a way to embed it in the permutation group.

By \newword{(combinatorial) unfolding} of a signed permutation $w$ of $\dbrack{n}$ we mean the following: to the right of $w$, attach a positionally reversed copy of $w$ with the signs flipped in every entry to get a permutation of $\pm \dbrack{n}$.   Conversely, given a permutation of $\pm\dbrack{n}$ where the $i$-th entry is the opposite sign of the $(2n+1-i)$-th entry, we can \newword{fold} the permutation to get a signed permutation on $\dbrack{n}$. For example, the unfolding of $[4,-2,3,-1]$ is 
$$[4,-2,3,-1,1,-3,2,-4].$$

If we normalize the entries of the unfolding, that is, if we replace $-\dbrack{n}$ with $\{n+1,\dots,2n\}$, then we get a map that assigns each signed permutation in $\dbrack{n}$ to a standard permutation in $S_{2n}$.

Given a (even) signed permutation $w$ in $B_n$ (or $D_n$), we let $\unfold(w)$ denote its unfolding $w'$ in $S_{2n}$. We let $\fold(w')$ be the reverse operation. Throughout the paper, we will slightly abuse the notation and identify a signed permutation of $\dbrack{n}$ with its unfolding in $\pm \dbrack{n}$. For signed permutations, the combinatorial unfolding is the same as the geometric one, whereas for even signed permutations, the combinatorial unfolding has nothing to do with the geometric version.

\begin{rmk}
   Note that the combinatorial unfolding of even signed permutations into permutations of $S_{2n}$ behaves quite strangely, especially compared to what happened for signed permutations in \cite{demAhopping}. For example, the relative ordering between the lengths of two elements is not preserved: we have $\ell(s_1) = \ell(s_n)$ in $D_n$ but if we take their combinatorial unfolding in $S_{2n}$ the former becomes a permutation that has length $2$, whereas the latter becomes a permutation that has length $4$.
\end{rmk}

The main technique used to treat the case of type $B_n$ in \cite{demAhopping} was the following lemma:

\begin{lem}[\cite{demAhopping}]
\label{lem:b_unfold}
For any signed permutations $w,v \in B_n$, we have 
\begin{equation}\label{eqn:fold-unfold}
w \star v = \fold(\unfold(w) \star \unfold(v)).
\end{equation}
\end{lem}

That is, we can analyze the Demazure product of signed permutations by studying the Demazure product of permutations obtained by unfolding them. This allows one to easily extend the results related to Demazure products for the case of type $A$ to type $B$. However, this is not the case for the type $D$ case. For example, look at the following example.

\begin{ex}
\label{ex:typeDbad}
Consider the (even) signed permutations $[1,4,-2,-3]$ and $[4,-1,2,-3]$. We have $[1,4,-2,-3] \star [4,-1,2,-3] = [-2,-1,4,3]$. However, if we take their unfoldings to obtain permutations $[1,4,-2,-3,3,2,-4,-1]$ and $[4,-1,2,-3,3,-2,1,-4]$ (upon normalizing, we can treat them as $[1,4,7,6,3,2,5,8]$ and $[4,8,2,6,3,7,1,5]$), and we take the Demazure product as elements of $S_8$, we get the permutation $[7,8,4,6,3,5,1,2]$. This folds back into $[-2,-1,4,-3]$, which does not even have even number of negative signs. 
\end{ex}

Fortunately, despite unfolding not behaving the way we want, the hopping operator can still emulate the Demazure product in the way we desire, and we will show that in the next section. 

\subsection{Maximal parabolic decomposition for type $D$}
\label{subsec:decomp}

In this subsection, we recall what the maximal parabolic quotient of type $D$ looks like. For $i \leq n-2$, let $S = \{s_i,\dots,s_n\}$ be a subset of generators of type $D_n$ (how the indicies are chosen is shown in \cref{fig:dn}). We get a parabolic subgroup $W_S$. Now, let $J = S \setminus \{s_i\}$. Then $W_J$ is a maximal parabolic subgroup of $W_S$ and we get $W^J_S$ as its right quotient. The Bruhat order on $W^J_S$ is drawn in \cref{fig:bruhatofquotient}.

\input{figures/typeDbruhatorder}


The elements of $W^J_S$ are one of the following four forms:
\begin{itemize}
    \item (form $0$) $id$,
    \item (form $1$) $s_i \dots s_j$ where $i \leq j \leq n-1$,
    \item (form $2$) $s_i \dots s_{n-2}s_{n}s_{n-1}\dots s_j$ where $i \leq j \leq n-1$,
    \item (form $3$) $s_i \dots s_{n-2}s_n$.
\end{itemize}

We also classify the elements of $W_{\{s_{n-1},s_n\}}$ using these forms:
\begin{itemize}
    \item (form $0$) $id$,
    \item (form $1$) $s_{n-1}$,
    \item (form $2$) $s_ns_{n-1}$,
    \item (form $3$) $s_n$.
\end{itemize}


Throughout the paper, we will classify the elements of $W_S^J$ and $W_{\{s_{n-1},s_n\}}$ using the above terminology. In the following example, we analyze how these elements act on the right of some $w \in D_n$ in one-line notation.

\begin{ex}
\label{ex:Dactright}
    Take $w = [2,-4,-1,5,3] \in D_5$. We will take a look at how $v$ acts on the right of $w$ and changes its one-line notation, depending on the form of $v$ classified as above. 
    \begin{itemize}
        \item (Form 1) if $v = s_i\dots s_j$ where $i \leq j \leq (n-1)$, then as we go from the one-line notation of $w$ to $wv$ we are cyclically shifting the entries in positions $i ,\dots, j+1$ to the left while preserving the signs. For example, $s_2s_3s_4$ cyclically shifts the entries in positions $2,3,4,5$, giving us
        $$ws_2s_3s_4 = [2,-1,5,3,-4].$$
        
        \item (Form 2) if $v = s_i \dots s_{n-2}s_ns_{n-1}\dots s_j$ where $i \leq n-2$ and $i \leq j \leq n-1$, then we are cyclically shifting the entries in positions $i,\dots, j$ to the left while preserving the signs, then flipping the signs of new entries in positions $j$ and $n$. For example, $s_2s_3s_5s_4$ shifts the entries in positions $2,3,4$, then flips the signs of new entries in positions $4$ and $5$, giving us
        $$ws_2s_3s_5s_4 = [2,-1,5,4,-3].$$ 
        If $v = s_{n}s_{n-1}$, then we simply flip the signs of the entries in positions $n-1$ and $n$.


        \item (Form 3) If $v = s_i\dots s_{n-2}s_n$, then we are cyclically shifting the entries in the positions $i,\dots, n$ to the left, then flipping the signs of the new entries in positions $n-1$ and $n$. For example, $s_2s_3s_5$ cyclically shifts the entries in positions $2,3,4,5$, then flips the signs of the new entries in positions $4$ and $5$, giving us
        $$ws_2s_3s_5 = [2,-1,5,-3,4].$$ 
        If $v = s_n$, then we cyclically shift the entries in positions $n-1,n$ to the left (or right, it is the same since there are only two entries), then flipping the signs of those entries.
        
    \end{itemize}
\end{ex}




We can repeatedly apply the above maximal parabolic decomposition by cutting off leaves in the order of $s_1,s_2,\dots,s_{n-2}$: first take $J_1 = S \setminus \{s_1\}$, then take $J_2 = J_1 \setminus \{s_2\}$ and so on until we are left with $W_{\{s_{n-1},s_{n}\}}$. In doing so, we get the following decomposition of any element in $D_n$.

\begin{prop}
\label{cor:decomp}
Given any element of $w \in D_n$, we can decompose $w$ into $Q_{n-1}Q_{n-2}\dots Q_1$ where $Q_i$ for each $i \leq n-2$ is an element of $W_{\{s_i,\dots,s_n\}}^{\{s_{i+1},\dots,s_n\}}$ and $Q_{n-1}$ is an element of $W_{\{s_{n-1},s_n\}}$.
\end{prop}

Consider $w = [2,-4,-1,5,3]$. Using \cref{cor:decomp}, we can decompose $w$ into $Q_4Q_3Q_2Q_1$ where we have $Q_4 = s_4s_5, Q_3 = s_3s_5, Q_2 = \text{id}, Q_1 = s_1s_2s_3s_4s_5s_3$. 




\input{sections/type_d.tex}

So we have ways to describe the Demazure product for type A,B,D Coxeter groups just using the one-line notation and not using reduced words. The natural next step would be to try this for affine permutations (type $\tilde{A}$).

\begin{question}
    Can one come up with hopping operators for affine permutations to describe the Demazure product?
\end{question}

\subsection*{Acknowledgments}
This research was carried out primarily during the 2024 Honors Summer Math Camp at Texas State University. The authors appreciate the support from the camp and also thank Texas State University for providing support and a great working environment. The authors also thank Edward Richmond for useful discussions.

\bibliography{biblio}
\bibliographystyle{siam}

\end{document}

%% file: settings.tex
\newif\ifdebug
\debugtrue
\usepackage[bookmarksnumbered,linktocpage,hypertexnames=false,colorlinks=true,linkcolor=blue,urlcolor=blue,citecolor=blue,anchorcolor=green,breaklinks=true,pdfusetitle]{hyperref}
\usepackage{graphicx,fullpage,xcolor}

\usepackage[utf8]{inputenc}        
\usepackage{lmodern}               
\usepackage[T1]{fontenc}           
\usepackage[english]{babel}

\usepackage{dsfont}
\usepackage[linesnumbered, ruled]{algorithm2e}

\usepackage{latexsym}              
\usepackage{amsmath}               
\usepackage{amssymb}               
\usepackage{amsfonts}              
\usepackage{amsthm}                

\usepackage{amsbsy}       
\allowdisplaybreaks
\usepackage{bbm}

\usepackage{tikz}                  
\usetikzlibrary{arrows,positioning,shapes}
\usepackage{tikzsymbols}

\usepackage{enumerate}             
\usepackage{enumitem}
\setlist{nolistsep}

\usepackage{xparse}                
\usepackage{hyperref}              
\usepackage{url}
\usepackage[capitalize,nameinlink]{cleveref}

\usepackage{color, xcolor}                 
\definecolor{darkgray}{rgb}{0.15,0.15,0.15}   
\definecolor{lightgray}{rgb}{0.94,0.94,0.94}  
\definecolor{lightlightgray}{rgb}{0.97,0.97,0.97}  
\definecolor{darkred}{rgb}{0.80,0.00,0.00}    
\definecolor{darkgreen}{rgb}{0.00,0.70,0.00}    
\definecolor{darkblue}{rgb}{0.00,0.00,0.70}    

\ifdebug
\fi

\usepackage{graphicx}              
\graphicspath{{fig/}}              
\usepackage{caption}               
\usepackage{subcaption}            
\usepackage{multicol}              
\usepackage{float}                 
\usepackage{framed}                
\usepackage[framemethod=tikz]{mdframed}
\usepackage{mathtools}
\usepackage{enumitem}              
\usepackage{glossaries}            

\usepackage[bottom]{footmisc}

\usepackage{tikz}          

\numberwithin{equation}{section}
\theoremstyle{plain}  
\newtheorem{thm}{Theorem}[section]
\newtheorem{lem}[thm]{Lemma}
\newtheorem*{lem*}{Lemma}
\newtheorem{prop}[thm]{Proposition}
\newtheorem*{prop*}{Proposition}

\newtheorem*{cor*}{Corollary}

\theoremstyle{definition}
\newtheorem{defn}[thm]{Definition}
\newtheorem*{defn*}{Definition}

\newtheorem*{notation*}{Notation}

\newtheorem*{problem*}{Problem}
\newtheorem{question}{Question}[section]
\newtheorem*{question*}{Question}

\theoremstyle{remark}  
\newtheorem{rmk}[thm]{Remark}
\newtheorem*{rmk*}{Remark}
\newtheorem{ex}[thm]{Example}
\newtheorem*{ex*}{Example}

\theoremstyle{plain}
\newtheorem{innercustomgeneric}{\customgenericname}
\providecommand{\customgenericname}{}
\newcommand{\newcustomtheorem}[2]{%
  \newenvironment{#1}[1]
  {%
   \renewcommand\customgenericname{#2}%
   \renewcommand\theinnercustomgeneric{##1}%
   \innercustomgeneric
  }
  {\endinnercustomgeneric}
}

\newcustomtheorem{customthm}{Theorem}
\newcustomtheorem{customlem}{Lemma}
\newcustomtheorem{customcor}{Corollary}

\providecommand*{\hr}[1][class-arg]{%
    \hspace*{\fill}\hrulefill\hspace*{\fill}
    \vskip 0.65\baselineskip
}

\NewDocumentCommand{\prob}{mg}{
    \IfValueTF{#2}{
        P(#1 \,\vert\, #2)
    }{
        P(#1)
    }
}

\newcommand\restr[2]{{
  \left.\kern-\nulldelimiterspace 
  #1 
  \right|_{#2} 
  }}

\newcommand{\Z}{\mathbb{Z}}

\newcommand{\dbrack}[1]{[[#1]]}
\newcommand{\newword}[1]{\textit{\textbf{#1}}}
\makeatletter
\newcommand\RedeclareMathOperator{%
  \@ifstar{\def\rmo@s{m}\rmo@redeclare}{\def\rmo@s{o}\rmo@redeclare}%
}
\newcommand\rmo@redeclare[2]{%
  \begingroup \escapechar\m@ne\xdef\@gtempa{{\string#1}}\endgroup
  \expandafter\@ifundefined\@gtempa
     {\@latex@error{\noexpand#1undefined}\@ehc}%
     \relax
  \expandafter\rmo@declmathop\rmo@s{#1}{#2}}
\newcommand\rmo@declmathop[3]{%
  \DeclareRobustCommand{#2}{\qopname\newmcodes@#1{#3}}%
}
\@onlypreamble\RedeclareMathOperator
\makeatother

\RedeclareMathOperator{\Re}{Re}

\DeclareMathOperator{\unfold}{unfold}
\DeclareMathOperator{\fold}{fold}

%% file: preamble.tex
\title{Demazure product and hopping in type D}

\author{Darren Han, Michelle Huang, Benjamin Keller, Suho Oh, Jerry Zhang}


\date{\today}
\maketitle

%% file: figures/typeDdiagram.tex
\begin{center}
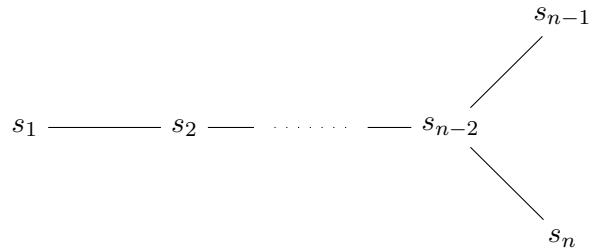

        \begin{tikzpicture}[scale = 0.7]
    \node (s1) at (-6,0) {$s_1$};
    \node (s2) at (-3,0) {$s_2$};
    \node (sn-2) at (2,0) {$s_{n-2}$};
    \node (sn-1) at (4.1,2.1) {$s_{n-1}$};
    \node (sn) at (4.1,-2.1) {$s_n$};
    \node (space1) at (-1.5, 0) {};
    \node (space2) at (.25, 0) {};
    \draw (s1) -- (s2);
    \draw (s2) -- (space1);
    \draw (space2) -- (sn-2);
    \draw (sn-2) -- (sn-1);
    \draw[loosely dotted] (space1) -- (space2);
    \draw (sn-2) -- (sn);
\end{tikzpicture}
    \captionsetup{width=1.0\linewidth}
  \captionof{figure}{The Coxeter diagram for $D_n$}
  \label{fig:dn}
\end{center}

%% file: figures/typeDbruhatorder.tex
\begin{center}
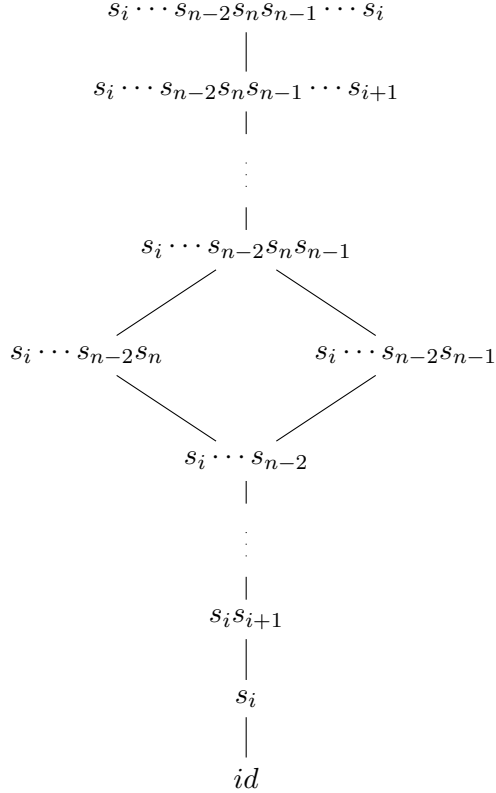

\begin{tikzpicture}[scale = 0.7]
    \node (id) at (0,-8) {$id$};
    \node (si) at (0,-6.5) {$s_i$};
    \node (sisi+1) at (0,-5) {$s_is_{i+1}$};
    \node (spaceb1) at (0, -4) {};
    \node (spaceb2) at (0, -3) {};
    \node (midbottom) at (0,-2) {$s_i\cdots s_{n-2}$};
    \node (midleft) at (-3,0) {$s_i\cdots s_{n-2}s_n$};
    \node (midright) at (3, 0) {$s_i\cdots s_{n-2}s_{n-1}$};
    \node (midtop) at (0, 2) {$s_i\cdots s_{n-2}s_ns_{n-1}$};
    \node (spacet1) at (0, 3) {};
    \node (spacet2) at (0, 4) {};
    \node (simidsi+1) at (0, 5) {$s_i\cdots s_{n-2}s_ns_{n-1}\cdots s_{i+1}$};
    \node (simidsi) at (0,6.5) {$s_i\cdots s_{n-2}s_ns_{n-1}\cdots s_i$};
    \draw (id) -- (si);
    \draw (si) -- (sisi+1);
    \draw (sisi+1) -- (spaceb1);
    \draw[loosely dotted] (spaceb1) -- (spaceb2);
    \draw (spaceb2) -- (midbottom);
    \draw (midbottom) -- (midleft);
    \draw (midbottom) -- (midright);
    \draw (midright) -- (midtop);
    \draw (midleft) -- (midtop);
    \draw (midtop) -- (spacet1);
    \draw[loosely dotted] (spacet1) -- (spacet2);
    \draw (spacet2) -- (simidsi+1);
    \draw (simidsi+1) -- (simidsi);
    
\end{tikzpicture}
    \captionsetup{width=1.0\linewidth}
  \captionof{figure}{Bruhat order of $W_{s_i,\ldots,s_n}^{s_{i+1},\ldots,s_n}$}
  \label{fig:bruhatofquotient}
\end{center}

%% file: sections/type_d.tex
\section{Hoppings in type $D$ and the main result}
In this section, we go over various properties of the hopping operator for even signed permutations, and state our main result: The Demazure product between two elements of Coxeter group $D_n$ can be expressed using only the one-line notation, by hopping on even signed permutations.

\subsection{Hopping in type $D$}
The hopping operator for type $D$ is actually defined in the exact same way as the hopping operator for type $B$ (since even signed permutations are signed permutations).


Beware that the hopping operator does not guarantee sending an even signed permutation to an even signed permutation. For example, consider $w = [1,5,-3,-4,2] \in D_n$. Let us try applying the hopping operator $h_{1,[3,-1,2]}$ to it:
$$h_{1, [3,-1,2]}([1, 5, -3, -4, 2]) = h_{1, [3,-1]}([2,5,-3,-4,1]) = h_{1,[3]}([2,5,-3,-4,-1]) = [2,5,-1,-4,-3].$$

The result $[2,5,-1,-4,-3]$ is not an even signed permutation, so it is no longer an element of $D_5$.

We label the positions of unfolding using $1,2,\dots,n,-n,\dots,-1$, following the ordering we placed on $\pm \dbrack{n}$. For any ordered list $L = [l_1,\dots,l_k]$ with entries from $\pm \dbrack{n}$ and an even signed permutation $w$, we define $w(L)$ to be the list $[w(l_1),\dots,w(l_k)]$ (with $w(-a) = -w(a)$).

We now investigate some properties of how the hopping operator interacts with the simple transpositions. For $s_1,\dots,s_{n-1}$, so everything except $s_n$, it interacts in a way that is exactly the same as what happened for the type $A$ and $B$ in \cite{demAhopping} as seen in the lemmas below. 

\begin{lem}[\cite{demAhopping}, Lemma 5]
\label{lem:hopirrele}
For $i \le n-2$ and $j > i$, we get
$$s_j h_{i,L}(s_{j}w) = h_{i, s_j(L)}(w).$$
\end{lem}
\begin{proof}
On the left-hand side, notice that since we apply $s_j$ left and right of the hopping operator, the roles of $j$ and $j+1$ are swapped. On the right-hand side, we are swapping the roles of $j$ and $j+1$ in the list used for the hopping operator. So we are looking at the same operation on both sides.
\end{proof}

\begin{lem}[\cite{demAhopping}, Lemma 6]
\label{lem:hoptrans}
For $i \le n-1$ and $L$ such that $i,i+1 \not \in L$, we have
$$s_ih_{i,L}(s_iw) = h_{i+1,s_i(L)}(w).$$ 
For $L$ such that $\pm (n-1), \pm n \not \in L$, we have
$$s_nh_{-(n-1),L}(s_nw) = h_{n,L}(w).$$

\end{lem}
\begin{proof}
  For the first claim, note that since we apply $s_i$ left and right of the hopping operator, the roles of $i$ and $i+1$ are swapped; all positions that get swapped in the process of hopping $i$ using $L$ are exactly the same as the positions that get swapped in the process of hopping $i+1$ using $s_i(L)$. For the second claim, the roles of $-(n-1)$ and $n$ are swapped.
\end{proof}


For example, let $w = [2,-4,-1,5,3] \in D_5$. Consider the actions of $s_2h_{2,[1,4,-5]}s_2$ and $h_{3,[1,4,-5]}$ on $w$. We have
$$[2,-4,-1,5,3] \xrightarrow[\quad s_2 \quad]{} [3,-4,-1,5,2] \xrightarrow[h_{2,[1,4,-5]}]{} [3,-2,-1,-4,-5] \xrightarrow[\quad s_2 \quad]{} [2,-3,-1,-4,-5],$$
and we also get
$$[2,-4,-1,5,3] \xrightarrow[h_{3,[1,4,-5]}]{} [2,-3,-1,-4,-5].$$
They yield the same result.

\subsection{Emulating the Demazure product with hopping}

In this section, we show how the Demazure product $v \star w$ can be emulated by the hopping operator when $v$ is simple: we start with the cases when $v$ is a simple transposition, then move onto cases where $v$ takes the form of elements classified in \cref{subsec:decomp}.

\begin{lem}
\label{lem:singletrans}
For $1 \le i \le n-1$, we get
\begin{equation}
\label{eq:form1fringe}
    s_i \star w = h_{i, [i+1]}(s_i w).
\end{equation}
We also have
    \begin{equation}
    \label{eq:form3fringe}
    s_n \star w = h_{n-1,[-n]} (s_nw).
    \end{equation}
\end{lem}
\begin{proof}
For the proof of \cref{eq:form1fringe}, notice that we only need to consider the relative order of $\pm (i),\pm (i+1)$. The only case where $s_i \star w$ and $s_i w$ differ (so $s_i \star w = w$) is when $i+1$ is left of $i$ inside the unfolding of $w$. When $i+1$ is on the right of $i$ within the unfolding of $w$, the operator $h_{i, [i+1]}$ is exactly the same as multiplying $s_i$ on the left.

For \cref{eq:form3fringe}, we only need to consider the relative order of $\pm (n-1), \pm n$, which we denote by $\chi$. Now look at the following table.

\begin{center}
\begin{tabular}{ |c|c|c| } 
  \hline $\chi$ & $s_n\star \chi$ & $h_{n-1,[-n]}(s_n\chi)$\\
  \hline\hline
  $n-1, n, -n, -(n-1)$  & $-n, -(n-1), n-1, -n$ & $-n, -(n-1), n-1, -n$\\
  \hline
  $n, n-1, -(n-1), -n$  & $-(n-1), -n, n, n-1$ & $-(n-1), -n, n, n-1$\\ 
  \hline
  $-(n-1), n, -n, n-1$  & $-(n-1), n, -n, n-1$ & $-(n-1), n, -n, n-1$\\
  \hline
  $n,-(n-1), n-1, -n$   & $-(n-1), n, -n, n-1$ & $-(n-1), n, -n, n-1$\\
  \hline
  $n-1, -n, n, -(n-1)$  & $-n, n-1, -(n-1), n$ & $-n, n-1, -(n-1), n$\\
  \hline
  $-n, n-1, -(n-1), n$  & $-n, n-1, -(n-1), n$ & $-n, n-1, -(n-1), n$\\ 
  \hline
  $-(n-1), -n, n, n-1$  & $-(n-1), -n, n, n-1$ & $-(n-1), -n, n, n-1$ \\
  \hline
  $-n, -(n-1), n-1, n$  & $-n, -(n-1), n-1, n$ & $-n, -(n-1), n-1, n$\\ 
  \hline
\end{tabular}
\end{center}
Since column $2$ and $3$ are the same, we get the desired claim.
\end{proof}

We now show how elements of form 1 interact with the hopping operator. The special case when $L=\emptyset$ in the following claim was done in \cite{demAhopping} for type A, but we have something that is stronger even for the type A case.

\begin{lem}
\label{cor:multtrans}
For any $w \in D_n$, $1 \leq i \leq j \leq n-1$ and $L$ such that $i,\dots,j+1 \not \in L$, we have
$$(s_i \dots s_j) \star h_{j+1,L}w = h_{i,[i+1,\dots,j+1,(s_i\dots s_j)(L)]} (s_i \dots s_j) w.$$
\end{lem}
\begin{proof}
We induct on $j-i$. For the base case of $i=j$, it follows from \cref{lem:singletrans} and \cref{eq:form1fringe}. Assume for the sake of induction that it holds for $i+1$. Then we have from the induction hypothesis and the base case that
\begin{align*}
(s_i \dots s_j) \star h_{j+1,L}w &= s_i \star h_{i+1,[i+2,\dots,j+1,(s_{i+1}\dots s_j)(L)]} (s_{i+1} \dots s_j)w \\
&= h_{i,[i+1,\dots,j+1,(s_i\dots s_j)(L)]} (s_i \dots s_j) w.
\end{align*}
This finishes the induction proof.
\end{proof}

We will create a particular list that will be used for the hopping operators, depending on the form of $Q_i$ (factors that come from the decomposition in \cref{cor:decomp}) we were looking at. We denote such a list as $L_i$ (beware that despite being indexed only using $i$, it actually depends on $Q_i$).
$$L_i := 
    \begin{cases}
        [i+1, \dots, j+1], &\text{if $Q_i$ is of form 1,} \\
        [i+1, \dots, n, -n, \dots, -(j+1)], &\text{if $i \leq n-2$ and $Q_i$ is of form 2,} \\
        [n,-n], &\text{if $i=n-1$ and $Q_{n-1}$ is of form 2,} \\
        [i+1, \dots, n-1, -n], &\text{if $i \leq n-2$ and $Q_i$ is of form 3,} \\
        [-n], &\text{if $i=n-1$ and $Q_{n-1}$ is of form 3.} \\
    \end{cases}
$$

\begin{prop}
\label{prop:hoppingtransfer}
For any even signed permutation $w$ have $$Q_i \star w = h_{i,L_i}(Q_iw).$$
\end{prop}
\begin{proof}
We have five cases, coming from the definition of $L_i$. The first case follows from \cref{cor:multtrans} by setting $L$ as an empty list. The fifth case is \cref{eq:form3fringe}. The third case follows from combining \cref{eq:form1fringe} and \cref{eq:form3fringe}. The fourth case follows from the fifth case by applying \cref{cor:multtrans} with $j=n-2$. So we only have to prove the second case.

We start with the claim that for $1 \le j \le i \le n-1$, we get
\begin{equation}
\label{eq:form1reverse}
s_i\dots s_j \star w = h_{-(i-1), [-(i-2), \dots, -j]}(s_i\dots s_j  w).  
\end{equation}
To show this, we induct on $i-j$. The base case for $j=i$ follows from \cref{eq:form1fringe} together with the observation that $h_{j,[j+1]}= h_{-(j+1), [-j]}.$ Using the same logic as the proof of \cref{lem:hoptrans}, we get
$$s_{j}h_{-j ,L}(s_jw) = h_{-(j+1),s_j(L)}(w),$$
as $s_j$ swaps the roles of $-j$ and $-(j+1)$. All that remains is to proceed with the induction in the exact same way as in the proof of \cref{cor:multtrans}, and the proof is omitted.

We then now aim to show the following.
\begin{equation}
\label{eq:n-1tojbase}
s_{n-1}s_ns_{n-2}\dots s_j\star w = h_{n-1,[n,-n,-(n-2),\dots,-j]}(s_{n-1}s_ns_{n-2}\dots s_jw).
\end{equation}

To show this, we do
\begin{align*}
\label{eq:n-1toj}
s_{n-1}s_ns_{n-2}\dots s_j\star w &= h_{n-1,[n,-n]}s_{n-1}s_nh_{-(n-1),[-(n-2),\dots,-j]}
(s_{n-2}\dots s_jw) \\ 
&= h_{n-1,[n,-n]} s_{n-1} h_{n,[-(n-2),\dots,-j]} (s_{n}s_{n-2}\dots s_jw) \\ 
&= h_{n-1,[n,-n]} h_{n-1,[-(n-2),\dots,-j]} (s_{n-1}s_n s_{n-2}\dots s_jw),
\end{align*}
where the first line follows from combining \cref{eq:form1reverse} and the third case, the second line follows from \cref{lem:hoptrans} and the third line follows from \cref{cor:multtrans}.

Finally, from \cref{eq:n-1tojbase} we can obtain the desired claim using \cref{cor:multtrans} with $j= n-2$.
\end{proof}



For example, letting $w = [2,-4,-1,5,3]$, consider $s_1s_2s_3s_5 \star w$. Then \cref{prop:hoppingtransfer} tells us that we can compute the Demazure product by taking $$h_{1,[2,3,4,-5]}(s_1s_2s_3s_5w) = h_{1,[2,3,4,-5]}([3,5,-2,-1,4]) = [3,-1,-2,5,4].$$



\subsection{Proof of the main result}



In this subsection, we present and prove our main result. We start by defining an operation that extracts a certain subword from the unfolding of an even signed permutation.

\begin{defn}
The lifting of $w \in D_n$ with respect to $i$, denoted by $w \nwarrow i$, is obtained by taking the subword strictly left of $i$ inside of the unfolding of $i$, but only take the entries between $(i,-i)$.
\end{defn}

For example, take $w = [2,-4,-1,5,3]$. Then $\unfold(w) = [2,-4,-1,5,3,-3,-5,1,4,-2]$. Then we have $w \nwarrow 4 = [5,-5]$ and $w \nwarrow 3 = [-4,5]$. 

Note that the definition of $w \nwarrow i$ for type $D$ is similar to that of type $B$, but this time we are taking entries in the interval $(i,-i)$ instead of $(i,-i]$. In the next lemma we show that for the hopping operator, if some $j$ and $-j$ appear right next to each other in the ordered list, we can freely swap them.

\begin{lem}
\label{lem:hopneg}
    For $i < j \leq n$, we have $$h_{i, [L, j, -j, L']} = h_{i, [L, -j, j, L']}.$$
\end{lem}
\begin{proof}
    It suffices to show that $h_{i,[j,-j]} = h_{i,[-j,j]}$. We verify this through a simple case-by-case analysis. We will only consider the relative order of $\pm i, \pm j$ in the unfolding as everything else is preserved, and we will denote this by $\chi$.
\begin{center}
    \begin{tabular}{|c|c|c|c|c|}
    \hline $\chi$ & $h_{i, [j,-j]}(\chi)$ & $h_{i, [-j, j]}(\chi)$\\
    \hline\hline $i,j,-j,-i$ &$-i,-j,j,i$ &$-i,-j,j,i$\\
    \hline $i,-j,j,-i$ & $-i,j,-j,i$ &$-i,j,-j,i$\\
    \hline $j,i,-i,-j$ &$-i,-j,j,i$ &$-i,-j,j,i$\\
    \hline $-j,i,-i,j$ &$-i,j,-j,i$ &$-i,j,-j,i$\\
    \hline $-i,j,-j,i$ &$-i,j,-j,i$ &$-i,j,-j,i$\\
    \hline $-i,-j,j,i$ &$-i,-j,j,i$ &$-i,-j,j,i$\\
    \hline $j,-i,i,-j$ &$-i,j,-j,i$ &$-i,j,-j,i$\\
    \hline $-j,-i,i,j$ &$-i,-j,j,i$ &$-i,-j,j,i$\\
    \hline
    \end{tabular}
\end{center}
From the table, we can conclude that the claim holds.
\end{proof}
For example, let $w = [2,-4,-1,5,3]$. Consider $h_{1, [3,2,-2,4]}w$ and $h_{1, [3,-2,2,4]}w$. In both cases, we get the same result as $[-1,-4,2,5,3]$.

Beware that swapping $j$ and $-j$ in a hopping operator can yield different results if $j$ and $-j$ are not right next to each other. For example, we have $h_{2, [4, -1, -4]}([1,2,3,4,5]) = [-2,4,3,-1,5]$ and $h_{2, [-4, -1, 4]}([1,2,3,4,5]) = [-2,-4,3,1,5]$ so we get different results.

\begin{defn}
    We define the following equivalence relation $\sim_i$ where $L \sim_i L'$ if $h_{i,L}(w) = h_{i,L'}(w)$ for all $w \in D_n$.
\end{defn}
For example, we may write \cref{lem:hopneg} as $[L,j,-j,L'] \sim_i [L,-j,j,L']$ whenever $i < j \leq n$. We now study some properties of the decomposition we defined in \cref{cor:decomp}. In the remainder of the paper we will use this decomposition when writing $w = Q_{n-1} \dots Q_1$. 

In general, when we compare $(Q_{n-1}\dots Q_{k+1}) \nwarrow i$ and $(Q_{n-1}\dots Q_k)\nwarrow i$, they are not necessarily the same. For example, when we take $w = [2,-4,-1,5,3]$, we get $Q_4 \nwarrow 4 = [-5,5] \neq [5,-5] = Q_4Q_3\nwarrow 4$.
However, we show that they are equivalent (behave in the same way when used in hopping operators).

\begin{lem}
\label{lem:addDk}
    For $k < i \leq (n-1)$, we have $$(Q_{n-1}\dots Q_{k+1}) \nwarrow i \sim_i (Q_{n-1}\dots Q_k)\nwarrow i.$$
\end{lem}
\begin{proof}
We will only consider how applying $Q_k$ on the right modifies the elements greater than $i$ in a one-line notation of some even signed permutation. Note that both $v := Q_{n-1}\dots Q_{k+1}$ and $v Q_k$ are the same as the identity permutation when restricted to entries $\{1,\dots,k-1\}$. Recall that in \cref{ex:Dactright} we gave a description of how $Q_k$ acts on the right of an even-signed permutation.

No matter what the form of $Q_k$ is, we first do a cyclic shift of the entries in $k,\dots,t$ (for some $t$, depending on the form) to the left. Since $k < i$, this process does not change the entries of $v \nwarrow i$ at all. Now for form $2$ or form $3$, we swap the signs of two entries (four in the unfolding), one of them being $k$ (which again is irrelevant in $v \nwarrow i$) and the other being $v(n)$. Since $v(n)$ and $v(-n)$ will be adjacent in $vQ_k \nwarrow i$, we get the desired result from \cref{lem:hopneg}.
\end{proof}

For example, let $w = [2,-4,-1,5,3]$, and $v = [4,3,-5,-1,2]$. We decompose $w = Q_4Q_3Q_2Q_1$ where we have $Q_4 = s_4s_5, Q_3 = s_3s_5, Q_2 = \text{id}$ and $Q_1 = s_1s_2s_3s_4s_5s_3$. We have $Q_4 \nwarrow 4 = [-5,5]$ and $Q_4Q_3 \nwarrow 4 = [5,-5]$.
For each of them of them used in a hopping operator we get:
$$h_{4, Q_4 \nwarrow 4}(v) = h_{4,[-5,5]}([4,3,-5,-1,2]) = [-4,3,5,-1,2],$$
$$h_{4, Q_4Q_3\nwarrow 4}(v) = h_{4,[5,-5]}([4,3,-5,-1,2]) = [-4,3,5,-1,2].$$
This aligns with \cref{lem:addDk} telling us that $Q_4 \nwarrow 4\sim_4 Q_4Q_3 \nwarrow 4$.

\begin{prop}
\label{prop:allformhop}
        
Take any $w = Q_{n-1} \dots Q_1 \in D_n$. For each $i$, let $L_i$ denote the list obtained from $Q_i$ as in \cref{prop:hoppingtransfer}. Then for any $i \leq n-2$, we have
    $$(Q_{n-1} \dots Q_{i+1})L_i\sim_i w\nwarrow i.$$
Also we have
$$ L_{n-1} \sim_{n-1} w \nwarrow (n-1).$$
    
\end{prop}
\begin{proof}
    For the first claim, it suffices to show $(Q_{n-1} \dots Q_{i+1})L_i \sim_i (Q_{n-1} \dots Q_i) \nwarrow i$ thanks to \cref{lem:addDk}. Notice that for $i \leq n-2$, we get $(Q_{n-1} \dots Q_{i+1})L_i = (Q_{n-1}\dots Q_i)T_i$ where $T_i$ is given by
$$T_i := 
    \begin{cases}
        [i, \dots, j], &\text{if $L_i = [i+1,\dots,j+1],$} \\
        [i+1, \dots,n-1, -n, n, -(n-1), \dots, -(j+1)], &\text{if $L_i = [i+1,\dots,n,-n,\dots,-(j+1)],$} \\
        [i+1, \dots, n-1], &\text{if $L_i = [i+1,\dots,n-1,-n]$.} \\
    \end{cases}
$$    
This can be explained by the following observation: what we had at positions $L_i$, after applying $Q_i$ on the right (recall from \cref{ex:Dactright} they all involve cyclic shift to the left on some positions, followed by a potential sign swap at two positions) becomes what we have at positions $T_i$. Also keep in mind that in $Q_{n-1} \dots Q_{i+1}$, the first $i$ entries are the same as the identity. So in each of the three cases, $i$ ends up exactly in position $j+1,-j,-n$ respectively in $Q_{n-1} \dots Q_i$. Hence, we find that $(Q_{n-1}\dots Q_i)T_i$ and $(Q_{n-1} \dots Q_i) \nwarrow i$ are exactly the same for the first and third cases (form 1 and form 3). For the second case (form 2), they are not necessarily the same (because we have $-n,n$ in $T_i$), but are equivalent due to \cref{lem:hopneg}.   

This finishes the proof of the first claim. For the second claim, we first do a simple casework check to compare $Q_{n-1} \nwarrow (n-1)$  and $L_{n-1}$:
\begin{center}
\begin{tabular}{|c|c|c|}
    \hline $Q_{n-1}$ & $Q_{n-1} \nwarrow (n-1)$ & $L_{n-1}$\\
    \hline\hline $\text{id}$ &$[]$ &$[]$\\
    \hline $s_{n-1}$ & $[n]$ &$[n]$\\
    \hline $s_ns_{n-1}$ &$[-n,n]$ &$[n,-n]$\\
    \hline $s_n$ &$[-n]$ &$[-n]$\\
    \hline
\end{tabular}
\end{center}

Together with \cref{lem:hopneg}, from the table we can conclude that $L_{n-1} \sim_{n-1}  Q_{n-1} \nwarrow (n-1)$. Since we have $w \nwarrow (n-1) \sim_{n-1} Q_{n-1} \nwarrow (n-1)$ from \cref{lem:addDk}, we are done with the proof of the second claim. 
\end{proof}
For example, let $w = [2,-4,-1,5,3]$, and $v = [-4,3,-5,-1,-2]$. We decompose $w = Q_4Q_3Q_2Q_1$ where we have $Q_4 = s_4s_5, Q_3 = s_3s_5, Q_2 = \text{id}$ and $Q_1 = s_1s_2s_3s_4s_5s_3$. We have $L_1 = [2,3,4,5,-5,-4]$ so $(Q_4Q_3Q_2)L_1 = [2,-4,-5,3,-3,5]$ and $w \nwarrow 1 = [2,-4,5,3,-3,-5]$. Now we consider $h_{1,[2,-4,-5,3,-3,5]}$ and $h_{1,[2,-4,5,3,-3,-5]}$ acting on $v$. We get
$$h_{1,[2,-4,-5,3,-3,5]}([-4,3,-5,-1,-2]) = [-4,-1,-5,3,-2],$$
$$h_{1,[2,-4,5,3,-3,-5]}([-4,3,-5,-1,-2]) = [-4,-1,-5,3,-2].$$
So we verify that $Q_4Q_3Q_2L_1 \sim_1 w\nwarrow 1$, as guaranteed by \cref{prop:allformhop}.

\begin{thm}[Demazure product in type D as hoppings]
\label{thm:mainD}
    For any $v,w \in D_n$,
    $$v\star w = h_{n-1, w \nwarrow n-1}\dots h_{1, w \nwarrow 1}(vw).$$
\end{thm}
\begin{proof}

From \cref{prop:hoppingtransfer}, we can rewrite 
\begin{equation}
\label{eq:usinghoptrans}
v\star w = h_{n-1, L_{n-1}}(Q_{n-1} h_{n-2, L_{n-2}}(Q_{n-2}\dots (Q_2 h_{1,L_1}(Q_1w)) ).
\end{equation}

Our goal is to transform the right-hand side of \cref{eq:usinghoptrans} into $h_{n-1, w \nwarrow n-1}\dots h_{1, w \nwarrow 1}(vw)$.


We get $Q_{n-1} \dots Q_{i+1} h_{i,L_i} = h_{i,Q_{n-1}\dots Q_{i+1}L_i}$ from \cref{lem:hopirrele}. And from \cref{prop:allformhop} we have $h_{n-1, w\nwarrow n-1} = h_{n-1, L_{n-1}}$ and $h_{i,w \nwarrow i} = h_{i,Q_{n-1}\dots Q_{i+1}L_i}$ for $i \leq n-2$. Then applying this on the left and moving inwards, we get
\begin{align*}
    v\star w &= h_{n-1, L_{n-1}}(Q_{n-1} h_{n-2, L_{n-2}}(Q_{n-2}\dots (Q_2 h_{1,L_1}(Q_1w)) ) \\ 
    &= h_{n-1,w\nwarrow n-1}(Q_{n-1} h_{n-2, L_{n-2}}(Q_{n-2}\dots (Q_2 h_{1,L_1}(Q_1w)) ) \\
    &= h_{n-1,w\nwarrow n-1}(h_{n-2, Q_{n-1}L_{n-2}}(Q_{n-1}Q_{n-2}\dots (Q_2 h_{1,L_1}(D_1w)) ) \\
    &= h_{n-1,w\nwarrow n-1}( h_{n-2, w\nwarrow n-2}(Q_{n-1}Q_{n-2}\dots (Q_2 h_{1,L_1}(Q_1w)) )\\
    &\qquad\qquad\vdots\\
    &= h_{n-1,w\nwarrow n-1}h_{n-2,w\nwarrow n-2}\dots h_{1,w\nwarrow 1} (vw).
\end{align*}
So we have transformed the right-hand side of \cref{eq:usinghoptrans} as we desire, and this finishes the proof.
\end{proof}

Notice that despite heavily relying on a specific decomposition using $Q_i$'s in the proof, our result does not use the decomposition nor the reduced words. For example, let $w = [2,-4,-1,5,3]$, and $v = [-4,3,-5,-1,-2]$. We get
\begin{align*}
h_{4, w\nwarrow 4}h_{3, w\nwarrow 3}h_{2, w\nwarrow 2}h_{1, w\nwarrow 1}(wv) 
&= h_{4, w\nwarrow 4}h_{3, w\nwarrow 3}h_{2, w\nwarrow 2}h_{1, w\nwarrow 1}([-5,-1,-3,-2,4]) \\
&= h_{4,[5,-5]}h_{3,[-4,5]}h_{1,[2,-4,5,3,-3,-5]}([-5,-1,-3,-2,4]) \\
&= h_{4,[5,-5]}h_{3,[-4,5]}([-1,-5,-3,-2,4]) \\
&= h_{4,[5,-5]}([-1,-3,-5,-2,4]) \\
&= [-1,-3,-4,-2,5].
\end{align*}

From \cref{thm:mainD}, we can conclude that 
$$w \star v = [-1,-3,-4,-2,5].$$